%% file: root.tex
\documentclass[journal,twoside,web,dvipsnames]{ieeecolor}
\usepackage{cite, amssymb, mathtools, amsfonts, graphicx, cuted,algorithm, algorithmic}
\usepackage{thmtools, thm-restate}
\usepackage{pgfplots, pgfplotstable}
    \usepgfplotslibrary{fillbetween, statistics}
    \pgfplotsset{compat=newest}

\newtheorem{proposition}{Proposition}
\newtheorem{assumption}{Assumption}
\newtheorem{definition}{Definition}

\def\logoname{LOGO-generic-web}
\definecolor{subsectioncolor}{rgb}{0,0.541,0.855}
\def\journalname{IEEE Transactions on Automatic Control}
\def\BibTeX{{\rm B\kern-.05em{\sc i\kern-.025em b}\kern-.08em
    T\kern-.1667em\lower.7ex\hbox{E}\kern-.125emX}}
\DeclareMathOperator{\vspan}{span}

\begin{document}
\title{Control Occupation Kernel Regression for Nonlinear Control-Affine Systems}

\author{Moad Abudia, Tejasvi Channagiri, Joel A. Rosenfeld, and Rushikesh Kamalapurkar 
    \thanks{This research was supported by the Air Force Office of Scientific Research (AFOSR) under contract numbers FA9550-20-1-0127 and FA9550-21-1-0134, and the National Science Foundation (NSF) under award numbers 2027976 and 2027999. Any opinions, findings and conclusions or recommendations expressed in this material are those of the author(s) and do not necessarily reflect the views of the sponsoring agencies.}
    \thanks{Moad Abudia is with the School of Mechanical and Aerospace Engineering, Oklahoma State University, Stillwater, OK 74078 USA (e-mail: abudia@okstate.edu).}
    \thanks{Tejasvi Channagiri is with the Department of Biostatistics, Harvard School of Public Health, Boston, MA 02115 USA (e-mail: tchannagiri@g.harvard.edu.).}
    \thanks{Joel A. Rosenfeld is with the Department of Mathematics and Statistics, University of South Florida, Tampa, FL 33620 USA (e-mail: tchannagiri@gmail.com,rosenfeldj@usf.edu).}
    \thanks{Rushikesh Kamalapurkar is with the Department of Mechanical and Aerospace Engineering, University of Florida, Gainesville, FL 32611 USA (e-mail: rkamalapurkar@ufl.edu).}%
}

\maketitle

\begin{abstract}
This manuscript presents an algorithm for obtaining an approximation of a nonlinear high order control affine dynamical system. Controlled trajectories of the system are leveraged as the central unit of information via embedding them in vector-valued reproducing kernel Hilbert space (vvRKHS). The trajectories are embedded as the so-called higher order control occupation kernels which represent an operator on the vvRKHS corresponding to iterated integration after multiplication by a given controller. The solution to the system identification problem is then the unique solution of an infinite dimensional regularized regression problem. The representer theorem is then used to express the solution as finite linear combination of these occupation kernels, which converts an infinite dimensional optimization problem to a finite dimensional optimization problem. The vector valued structure of the Hilbert space allows for simultaneous approximation of the drift and control effectiveness components of the control affine system. Several experiments are performed to demonstrate the effectiveness of the developed approach.
\end{abstract}

\section{Introduction}\label{sec:introduction}

 Consider a dynamical system of the form $\dot x = f(x)$ with unknown dynamics, $f: \mathbb{R}^n \to \mathbb{R}^n$. Given an observed trajectory, $\gamma:[0,T] \to \mathbb{R}^n$ satisfying the dynamics, system identification routines for dynamical systems have traditionally relied on taking numerical derivatives of observed trajectory data \cite{SCC.Chowdhary.Johnson2011,SCC.Brunton.Proctor.ea2016a}, or numerical gradient of a cost function 
  \cite{SCC.Ljung1998}. However, numerical derivatives are sensitive to signal noise, where an addition of white noise can cause unbounded disturbances for numerical differentiation. Methods such as the SINDy \begingroup \color{red} \cite{SCC.Brunton.Proctor.ea2016a}  and SINDYc algorithm \cite{SCC.Brunton.Proctor.ea2016} \endgroup ameliorate this issue through a total variation regularization, but even these methods have limitations with the respect to noise.

Kernel methods developed by the machine learning community have been adopted for system identification purposes by the control community \cite{pillonetto2014kernel,pillonetto2010new,pillonetto2011new}. The complexity of kernel methods is further investigated in the system identification context in  \cite{pillonetto2015tuning}, where it is shown that the tuning of hyperparameters is a problem that still persists today.

Recently, a collection of results surrounding the concept of occupation kernels have appeared, where system identification problems are addressed not through numerical differentiation, but rather through integration \cite{SCC.Rosenfeld.Russo.ea2024}. While the analysis in this paper does not explicitly consider measurement noise, the developed algorithm relies on numerical integration, and as indicated by the simulation results, the integration approach is less sensitive to signal noise compared to state-of-the-art methods such as SINDYc. The integration approach can be incorporated in system identification routines naturally through reproducing kernel Hilbert spaces (RKHSs) and occupation kernels \cite{SCC.Rosenfeld.Russo.ea2024}. 

Occupation kernels are functions in a RKHS that represent trajectories of a dynamical system. Since the span of the occupation kernels can be shown to be dense in the RKHS, they can be leveraged as basis functions for approximation (see Proposition \ref{prop:SystemTrajectoryDensity}). Indeed, occupation kernels have been used in precisely that manner for motion tomography in \cite{SCC.Russo.Kamalapurkar.ea2022} as well as for a regression approach to fractional order nonlinear system identification in \cite{SCC.Li.Rosenfeld2020}. Occupation kernels are also leveraged as basis functions for the construction of eigenfunctions for finite rank representations of Liouville operators in a continuous time dynamic mode decomposition routine in \cite{SCC.Rosenfeld.Kamalapurkar.ea2022}. Occupation kernels generalize the idea of occupation measures, which have been leveraged extensively in optimal control, in a way analogous to how kernel functions generalize delta distributions \cite{SCC.Lasserre.Henrion.ea2008,SCC.Majumdar.Vasudevan.ea2014,SCC.Korda.Henrion.ea2016}. Like kernel functions, occupation kernels provide a function theoretic analog of their respective measures in a RKHS.

The present manuscript generalizes two different approaches to system identification using occupation kernels, occupation kernels for higher order control affine systems are introduced in Section \ref{sec:kernels}, and are leveraged as basis functions for the resolution of a regularized regression problem for higher order control affine systems of the form $\frac{\mathrm{d}^{s}}{\mathrm{d}t^{s}} x = f(x) + g(x) u$. Specifically, this paper extends the regression approach for fractional order dynamical systems in \cite{SCC.Li.Rosenfeld2020} as well as giving a generalization of the control occupation kernels used in combination with control Liouville operators in \cite{SCC.Rosenfeld.Kamalapurkar2021}.

System identification using control occupation kernel regression (COKR) is realized through a regularized regression formulation, where the regularization forces the minimizer to be a linear combination of the occupation kernels, and in turn, leads to a finite dimensional optimization problem resolved by a matrix equation given in Section \ref{sec:regression-for-control}. In Section \ref{sec:FRR}, it is shown that COKR can be viewed as a generalization of kernel ridge regression (KRR) (see, for example, \cite[Appendix A]{zhdanov2013identity}) and that the convergence properties of KRR in terms of function evaluation can be extended to COKR as convergence properties in terms of inner product evaluation (see Proposition \ref{prop:prediction loss}).  


\section{Problem Statement}
\label{sec:problem}

The objective of this manuscript is to learn an unknown higher order control affine system
\begin{equation}\label{eq:control-affine-dynamics}
\frac{\mathrm{d}^{s}x}{\mathrm{d}t^{s}} = f(x) + g(x) u
\end{equation}
from observed piecewise continuous control signals $\{ u_j : [0,T_j] \to \mathbb{R}^m \}_{j=1}^M$ and corresponding Carath\'{e}odory solutions \cite[Section 3.1]{khalil2002nonlinear} $\{ \gamma_{u_j}:[0,T_j] \to \mathbb{R}^n\}_{j=1}^M$ of \eqref{eq:control-affine-dynamics}, where $f : \mathbb{R}^n \to \mathbb{R}^n$ is the drift function, $g:\mathbb{R}^n \to \mathbb{R}^{n\times m}$ is the control effectiveness matrix, and $s \in \mathbb{N}$ is the order of the system, which is assumed to be known.

For several selections of $f$, $g$, and $s$, the formulation in \eqref{eq:control-affine-dynamics} specializes to several different system identification problems. When $g \equiv 0$, the problem reduces to determining the unknown dynamics $f$ for a higher order dynamical system. When $s = 1$ the problem becomes a system identification problem for a first order control-affine system. If $s=1$ and $g \equiv 0$, the method developed in this manuscript agrees with that of \cite{SCC.Li.Rosenfeld2020} for integer-order systems with $q=1$.

Systems of the form \eqref{eq:control-affine-dynamics} encompass $s-$th order linear systems and Euler-Lagrange models with invertible inertia matrices, and hence, represent a wide class of physical plants, including but not limited to robotic manipulators and autonomous ground, aerial, and underwater vehicles.

In order to facilitate the description of the controlled dynamical system in terms of operators, a vector valued Reproducing Kernel Hilbert Space (vvRKHS) framework is utilized in this paper. 

\section{Higher Order Occupation Kernels}
\label{sec:kernels}

Given a Hilbert space $\mathcal{Y}$ and a set $X$, a vector valued RKHS, $H$, is a Hilbert space of functions from $X$ to $\mathcal{Y}$, such that for each $v\in \mathcal{Y}$ and $x \in X$, the functional that maps $f \in H $ to $\langle f(x), v \rangle_{\mathcal{Y}} \in \mathbb{R}$ is bounded. Hence for each $x \in X$ and $v \in \mathcal{Y}$, there is a function $K_{x, v} \in H$ such that $\langle f(x), v \rangle_{\mathcal{Y}} = \langle f, K_{x,v} \rangle_{H}$. The mapping $v \mapsto K_{x,v}$ is linear over $\mathcal{Y}$; hence, $K_{x,v}$ can be expressed as an operator over $\mathcal{Y}$ \begingroup\color{red} as $K_x v := K_{x,v}$, where $K_x:\mathcal{Y}\to H$ is called the kernel operator of $H$ centered at $x$\endgroup. The operator $K(x,y) := K_y^*K_x$ is called the reproducing kernel operator of $H$, where $K_y^*$ is the adjoint of $K_y$.

 In this paper, unless otherwise specified, we assume that $\mathcal{Y} = \mathbb{R}^{m+1}$ (viewed as row vectors), $X=\mathbb{R}^n$, and $H$ is a vvRKHS of $\mathcal{Y}-$valued continuous functions over $X$. Since $\mathcal{Y}$ is the space of row-vectors, the operation on $v \in \mathbb{R}^{m+1}$ by the kernel operator $K_x$ will be expressed as $v K_{x}$. A trivial extension of the proof of \cite[Proposition 2]{rosenfeld2024dynamic} shows that \begingroup \color{red} for any vvRKHS $H$ of continuous functions, \endgroup given any continuous signal, $\theta : [0,T] \to \mathbb{R}^n$ and any bounded measurable signal $u : [0,T] \to \mathbb{R}^m$, the functional $\mathcal{T}_{\theta ,u}^s:H\to\mathbb{R}$, defined as $\mathcal{T}_{\theta ,u}^s h = \frac{1}{(s-1)!}\int_0^T (T-t)^{s-1} h(\theta(t))\begin{pmatrix} 1 \\ u(t) \end{pmatrix}\mathrm{d}t$ for all $h\in H$ is bounded. Control occupation kernels were introduced in \cite{SCC.Rosenfeld.Kamalapurkar2021} as the unique functions in $H$ that represent the functional $\mathcal{T}_{\theta ,u}^1$, and can be generalized to higher order dynamical systems through the Cauchy iterated integral formula. Note that Definition \ref{def:control_occupation_kernel} and Proposition \ref{prop:evaluate-occupation} require weaker regularity assumptions than those made in the problem statement, which ensure existence and uniqueness of solutions of \eqref{eq:control-affine-dynamics}. 
\begin{definition}\label{def:control_occupation_kernel}
Given a continuous signal $\theta:[0,T]\to\mathbb{R}^n$ and a bounded measurable signal $u:[0,T]\to\mathbb{R}^m$, the \emph{control occupation kernel of order $s \in \mathbb{N}$ corresponding to $\theta$} in $H$, denoted by  $\Gamma_{\theta,u}^{(s)}$ is given as the unique function that represents the bounded functional $\mathcal{T}_{\theta ,u}^s$ as
$\langle h, \Gamma_{\theta,u}^{(s)} \rangle_H = \mathcal{T}_{\theta ,u}^s h.$
\end{definition}
This definition enables the treatment of particular higher order dynamical systems without state augmentation, i.e., without writing the system in terms of first order ordinary differential equations. State augmentation for systems studied in this note would result in needless approximation of already known parts of the dynamics.

\begin{proposition}\label{prop:evaluate-occupation}
Given a continuous signal $\theta:[0,T]\to\mathbb{R}^n$, a bounded measurable signal $u:[0,T]\to\mathbb{R}^m$, and $s \in \mathbb{N}$, \begin{equation}\label{eq:evaluate-occupation}
    \Gamma^{(s)}_{\theta,u}(x) = \frac{1}{(s-1)!} \int_0^T (T-t)^{s-1} \begin{pmatrix} 1 & u^{\top}(t) \end{pmatrix} K_{\theta(t)}(x)\mathrm{d}t,
\end{equation}
\begingroup\color{red} where $K_{\theta(t)}$ is the kernel operator of $H$ centered at $\theta(t)$.\endgroup
\end{proposition}

\begin{proof}
Consider $\langle \Gamma^{(s)}_{\theta,u}(x), v \rangle_{\mathcal{Y}} = \langle \Gamma^{(s)}_{\theta,u},vK_{x} \rangle_H,$ which follows from the definition of a vvRKHS. Symmetrically, since $H$ is a real-valued vvRKHS, it follows that
\begingroup\medmuskip=0mu\thinmuskip=0mu\thickmuskip=0mu
\begin{gather*} 
    \langle \Gamma^{(s)}_{\theta,u}(x), v \rangle_{\mathcal{Y}} = \langle \Gamma^{(s)}_{\theta,u},vK_{x} \rangle_H = \langle vK_{x}, \Gamma^{(s)}_{\theta,u} \rangle_H =\\ \frac{1}{(s-1)!} \int_0^T (T-t)^{s-1} [v K_{x}](\theta(t)) \begin{pmatrix} 1 \\ u(t) \end{pmatrix} \mathrm{d}t\\
    = \left\langle \frac{1}{(s-1)!} \int_0^T (T-t)^{s-1} \begin{pmatrix} 1 & u(t)^{\top} \end{pmatrix} K_{\theta(t)}(x)\mathrm{d}t, v \right\rangle_{\mathbb{R}^{m+1}}.
\end{gather*}
\endgroup
Hence, the equality holds for all $v \in \mathbb{R}^{m+1},$ and the result follows.
\end{proof}


Due to the presence of an integral in \eqref{eq:evaluate-occupation}, a numerical integration routine, such as Simpson's rule, is needed to evaluate the occupation kernels and inner products involving the occupation kernels.

The following proposition shows that Definition \ref{def:control_occupation_kernel} can be used to compute inner products in the vvRKHS. The notation $(\cdot)_j$ is used to denote the $j-$th row of the matrix $(\cdot)$.

\begin{proposition}\label{prop:ftc}
If $\gamma_u: [0,T] \to \mathbb{R}^n$ is a Carath\'{e}odory solution of \eqref{eq:control-affine-dynamics} with $\begin{pmatrix} (f)_{j} & (g)_{j} \end{pmatrix} \in H$ for each $j=1,\ldots,n$, under a piecewise continuous controller $u:[0,T]\to\mathbb{R}^m$, starting from the initial conditions $\left\{\frac{\mathrm{d}^{\ell}\gamma_{u}}{\mathrm{d}t^{\ell}}(0)\right\}_{\ell = 0}^{s-1}$, then for all $s\in\mathbb{N}$,
$
    \langle \begin{pmatrix} (f)_{j} & (g)_{j} \end{pmatrix}, \Gamma^{(s)}_{\gamma_u, u} \rangle_H = \left(\gamma_{u}(T) - \sum_{\ell = 0}^{s-1} \frac{T^\ell}{\ell!}\frac{\mathrm{d}^{\ell}\gamma_{u}}{\mathrm{d}t^{\ell}}(0)\right)_j.
$
\end{proposition}

\begin{proof}
Since $\begin{pmatrix} (f(\gamma_u(t)))_{j} & (g(\gamma_u(t)))_{j} \end{pmatrix} \begin{pmatrix}1 \\ u(t)\end{pmatrix} = (f(\gamma_u(t)))_{j} + (g(\gamma_u(t)))_{j} u(t) =  \left(\frac{\mathrm{d}^s\gamma_u}{\mathrm{d}t^s}(t)\right)_j$ for almost all $t\in [0,T]$, and $\Gamma_{\gamma_u,u}^{(s)}$ implements Cauchy's iterated integral formula through the inner product of the Hilbert space, the proposition follows through iterated application of the fundamental theorem of calculus.
\end{proof}


Section \ref{sec:regression-for-control} explores a regularized regression approach for learning higher order control affine dynamical systems, where the dynamics can be expressed as a linear combination of higher order control occupation kernels. Therefore, it is necessary to have a cogent method for evaluating the higher order control occupation kernels, and this can be realized through inner products with the kernels $e_j K_{x}$, where $e_j$ is the $j$-th cardinal basis function for $\mathbb{R}^{m+1}$. The following proposition provides an explicit formulation to compute the control occupation kernel in \eqref{eq:evaluate-occupation}.

If variable length trajectories are admitted, then
it can be shown that the span of the occupation kernels corresponding to trajectories $\theta_u$ that result from the application of the control $u$ to the system in \eqref{eq:control-affine-dynamics} is dense in $H$.
\begin{restatable}[]{proposition}{densityProposition}\label{prop:SystemTrajectoryDensity}
For any order $s\in\mathbb{N}$ and any $\{\gamma_\ell\}_{\ell=1}^{s-1}\subset\mathbb{R}^n$, the span of the set
\[
    A_s\coloneqq\left\{\Gamma_{\gamma_{u,\gamma_0},u}^{(s)}\mid \begin{gathered}
        \gamma_0\in\mathbb{R}^n,T_u\in[0,T],\\
        u\in \mathcal{PC}^0([0,T_u];\mathbb R^m)
    \end{gathered} \right\}
\]
is dense in $H$, where $\gamma_{u,\gamma_0}$ is a Carath\'{e}odory solution of \eqref{eq:control-affine-dynamics} under the control $u$ and starting from the initial conditions $\frac{\mathrm{d}^j}{\mathrm{d}t^j}\gamma_{u,\gamma_0}(0) = \gamma_j$ for $j=0,\ldots,s-1$, \begingroup\color{red} $\Gamma_{\gamma_{u,\gamma_0},u}^{(s)}$ denotes the control occupation kernel of order $s$ corresponding to the signal $\gamma_{u,\gamma_0}$ in $H$, \endgroup and $\mathcal{PC}^0$(A;B) denotes the space of piecewise continuous functions with domain $A$ and co-domain $B$.
\end{restatable}

\begin{proof}
    See the appendix.
\end{proof}
Proposition \ref{prop:SystemTrajectoryDensity} motivates the use of control occupation kernels for system identification. If $H$ is a vvRKHS of a universal kernel \cite[Definition 2]{carmeli2010vector}, then any continuous function can be approximated, uniformly over any compact set, by a function in $H$ \cite[Theorem 3]{carmeli2010vector}, and any function in $H$ can be approximated using a linear combination of sufficiently many control occupation kernels.

\section{A Regression Method using Control Occupation Kernels for Control Affine Systems}
\label{sec:regression-for-control}
The COKR method developed in this manuscript solves a regularized regression problem posed over the vvRKHS, $H$. Through the representer theorem, it is shown that that higher order control occupation kernels arise naturally as a collection of basis functions for approximating the control affine dynamics. 
\begin{multline*}
    \left\langle \begin{pmatrix} (f)_{i} & (g)_{i} \end{pmatrix}, \Gamma^{(s)}_{\gamma_u, u} \right\rangle_H = \left(\gamma_{u}(T) - \sum_{\ell = 0}^{s-1} \frac{T^\ell}{\ell!}\frac{\mathrm{d}^{\ell}\gamma_{u}}{\mathrm{d}t^{\ell}}(0)\right)_i\\
    =\frac{1}{(s-1)!}\int_0^{T} (T - t)^{s-1} \left( (f(\gamma_{u}(t)))_{i} + (g(\gamma_{u}(t)))_{i} u(t)\right) \mathrm{d}t,
\end{multline*}
given $\lambda > 0$ and the controllers and controlled trajectories from Section \ref{sec:problem}, the regularized regression problem to determine an approximation of the $i$-th row of $f$ and $g$ within the vector valued RKHS is given as
\begin{gather}
\min_{\begin{pmatrix} (\hat f)_{i} &  (\hat g)_{i}\end{pmatrix} \in H}  \sum_{j=1}^M \Bigg[ \left\langle \begin{pmatrix}  (\hat f)_{i} &  (\hat g)_{i} \end{pmatrix}, \Gamma^{(s)}_{\gamma_{u_j}, u_j} \right\rangle_H  \nonumber \\
-\left( \gamma_{u_j}(T_j) - 
\sum_{\ell=1}^{s-1}\frac{T_j^\ell}{\ell!}\frac{\mathrm{d}^{\ell}\gamma_{u_j}}{\mathrm{d}t^{\ell}}(0)\right)_i\Bigg]^2 + \lambda \| \begin{pmatrix} (\hat f)_{i} &  (\hat g)_{i}\end{pmatrix} \|^2_H,\label{eq:regression}
\end{gather}
where $\lambda > 0$ is a regularization parameter, and $
    \frac{1}{(s-1)!}\int_0^{T_j} (T_j - t)^{s-1} \left((\hat{f})_{i}(\gamma_{u_j}(t)) +  (\hat g)_{i}(\gamma_{u_j}(t)) u_j(t)\right) \mathrm{d}t$
is equal to the inner product $\left\langle  \begin{pmatrix} (\hat f)_{i} &  (\hat g)_{i}\end{pmatrix}, \Gamma_{\gamma_{u_j},u_j}^{(s)} \right\rangle_H$. Using the Representer Theorem for occupation kernels (cf. \cite{SCC.Kimeldorf.Wahba1970} and  \cite[Proposition 1]{SCC.Li.Rosenfeld2020}), the minimizer of \eqref{eq:regression} can be expressed as a linear combination of occupation kernels $\begin{pmatrix} (\hat f)_{i} &  (\hat g)_{i}\end{pmatrix} = w_{i,1} \Gamma_{\gamma_{u_1},u_1}^{(s)}+ \cdots + w_{i,M} \Gamma^{(s)}_{\gamma_{u_M},u_M}.
$
Using this representation of the minimizer, the inner products in the optimization problem can be computed as$
    \langle \begin{pmatrix} (\hat f)_{i} &  (\hat g)_{i}\end{pmatrix}, \Gamma_{\gamma_{u_j},u_j}^{(s)} \rangle_H \\=\left\langle \sum_{k=1}^M w_{i,k} \Gamma_{\gamma_{u_k},u_k}^{(s)}, \Gamma_{\gamma_{u_j},u_j}^{(s)} \right\rangle_H=
    (w)_i(G)^{j}$

where $(w)_i \coloneqq \begin{pmatrix} w_{i,1}&\ldots&w_{i,M} \end{pmatrix}$ denotes the $i-$the row of the weight matrix $w\in\mathbb{R}^{n\times M}$, $(G)^{j}$ denotes the $j-$th column of the Gram matrix $G\in\mathbb{R}^{M\times M}$, whose $i,j-$th element, denoted by $G_{i,j}$, is given by $G_{i,j} = \left\langle\Gamma^{(s)}_{\gamma_{u_i},u_i},\Gamma^{(s)}_{\gamma_{u_j},u_j}\right\rangle_H $, where 
\begin{multline}
    \left\langle \Gamma^{(s)}_{\gamma_{u_i},u_i}, \Gamma^{(s)}_{\gamma_{u_j},u_j} \right\rangle_{H}
    =\left( \frac{1}{(s-1)!}\right)^2  \int\limits_{0}^{T_j}\int\limits_{0}^{T_i} (T_i-t)^{s-1} \\(T_j-t)^{s-1} \begin{bmatrix}1 & u_i^{\top}(\tau)\end{bmatrix}  K\left(\gamma_{u_j}(t),\gamma_{u_i}(\tau)\right)\begin{bmatrix}
    1 \\ u_j(t)
    \end{bmatrix}\mathrm{d}\tau\, \mathrm{d}t.\label{eq:Control_occ_ker_gram_matrix}
\end{multline}
Furthermore, the norm in the optimization problem can be computed as $\| \begin{pmatrix} (\hat f)_{i} &  (\hat g)_{i}\end{pmatrix} \|^2_H = (w)_i G (w)_i^\top$. Hence, the resolution \eqref{eq:regression} reduces to the finite dimensional convex optimization problem
\begin{equation} 
    \min_{(w)_i \in \mathbb{R}^{1\times M}} \sum_{j=1}^M\left( \left[ (w)_i(G)^{j} - (D)_{i}^{j}\right]^2 + \lambda (w)_i G (w)_i^\top\right). \label{eq:RegularizedRegressionFinite}
\end{equation}
where $(D)_{i}^{j}$ is the $i,j-$th element of the end point difference matrix $D\in\mathbb{R}^{n\times M}$, whose $j-$th column is $(D)^j = \gamma_{u_j}(T_j) - \sum_{\ell=1}^{s-1}\frac{T_j^\ell}{\ell!}\frac{\mathrm{d}^{\ell}\gamma_{u_j}}{\mathrm{d}t^{\ell}}(0)$.

The optimization problem in \eqref{eq:RegularizedRegressionFinite} is written with respect to a single dimension in the state space. Once all the state space dimensions are concatenated, the combined problem can be solved via the resolution of the linear system $\left(G+\lambda I_M\right)w^\top=D^\top$.

The resultant approximation is given as
$\begin{pmatrix} \hat{f}(x) & \hat{g}(x) \end{pmatrix} = \sum_{j=1}^M (w)^j \Gamma^{(s)}_{\gamma_{u_j},u_j}(x)$
where $(w)^j$ denotes the $j-$th column of $w$. The COKR technique is summarized in Algorithm \ref{alg:COKR}.
 
\begin{algorithm}\caption{\label{alg:COKR}The COKR algorithm}\begin{algorithmic}[1]
        \renewcommand{\algorithmicrequire}{\textbf{Input:}}
        \renewcommand{\algorithmicensure}{\textbf{Output:}}
      
        \REQUIRE Trajectories $\{\gamma_{u_i}\}_{i=1}^{M}$, control signals $\{u_i\}_{i=1}^{M}$, a regularization parameter $\lambda$, the order of the system $s$, a numerical integration procedure, and a reproducing  kernels  $K$ of  $H$.
        \ENSURE $\begin{pmatrix} \hat{f}(x) & \hat{g}(x) \end{pmatrix} $
        \STATE $G \leftarrow \left(\left\langle \Gamma_{\gamma_{u_i},u_i}^{(s)}, \Gamma_{\gamma_{u_j},u_j}^{(s)} \right\rangle_{H} \right)_{i,j=1}^M$; using \eqref{eq:Control_occ_ker_gram_matrix} 
        \STATE $ D \leftarrow \left( \left(\gamma_{u_j}(T_j)\right)_i - \left(\sum_{\ell=1}^{s-1}\frac{T_j^\ell}{\ell!}\frac{\mathrm{d}^{\ell}\gamma_{u_j}}{\mathrm{d}t^{\ell}}(0)\right)_i\right)_{i,j=1}^{n,M}$
        \STATE $w\leftarrow ((G+\lambda I_M)^{-1} D^\top)^\top$ 
        \STATE $\begin{pmatrix} \hat{f}(x) & \hat{g}(x) \end{pmatrix} \leftarrow \sum_{j=1}^M (w)^j \Gamma^{(s)}_{\gamma_{u_j},u_j}(x)$; using \eqref{eq:evaluate-occupation}
        \RETURN $\begin{pmatrix} \hat{f}(x) & \hat{g}(x) \end{pmatrix}$ 
    \end{algorithmic} 
\end{algorithm}

\section{Functional ridge regression} \label{sec:FRR}

In the developed COKR method, we are given a set of functions $\{\Gamma_{\gamma_{u_i},u_i}^{s}\}_{i=1}^M$ and a set of constants $(D)_{i}^{j}$ that satisfy $\left\langle (h)_i, \Gamma^{(s)}_{\gamma_{u_j}, u_j} \right\rangle_H = (D)_{i}^{j}$ with $(h)_i = \begin{pmatrix} (f)_{i} &  (g)_{i}\end{pmatrix} $. We then proceed to solve \begingroup\small$
    \min_{(h)_i \in H}  \sum_{j=1}^M \left(\left[ \left\langle (h)_i, \Gamma^{(s)}_{\gamma_{u_j}, u_j} \right\rangle_H - (D)_{i}^{j}\right]^2 + \lambda \| (h)_i \|^2_H\right).$\endgroup

Using the representer theorem, we can also show that given a new function $\Gamma_{\gamma_{u},u}^{s}$, the COKR predictor $(\hat{h})_i$ satisfies $(\delta)_i \coloneqq\langle (\hat{h})_{i},\Gamma_{\gamma_{u},u}^{s}\rangle_{H} = (D)_i(G + \lambda I_M)^{-1}\begin{pmatrix} \langle \Gamma^{(s)}_{\gamma_{u_1}, u_1}, \Gamma^{(s)}_{\gamma_{u}, u} \rangle_{H} & \ldots & \langle \Gamma^{(s)}_{\gamma_{u_M}, u_M}, \Gamma^{(s)}_{\gamma_{u}, u} \rangle_{H} \end{pmatrix}^{\top}$. The COKR method can therefore be seen as a generalization of kernel ridge regression (KRR) (see, for example, \cite[Appendix A]{zhdanov2013identity}), where the control occupation kernels replace the reproducing kernels. 

The developed generalization is useful in problems where direct samples of functions to be approximated are not available for training. If inner products of the unknown functions against a set of basis functions can be computed, and if that set of basis functions satisfies the representer theorem, then the generalized KRR method can be used to generate useful estimates of the unknown functions. The affine system identification problem studied in this paper is one example of approximation problems that have this property. The price paid for the generalization is that instead of performance guarantees in terms of function evaluation, we  get performance guarantees in terms of inner product evaluation against the set of basis used for the representation. Performance in terms of function evaluation then depends on the selected set of basis, and needs to be analyzed separately.

Performance guarantees developed for KRR in terms of function evaluation can be interpreted as performance guarantees in terms of inner product evaluation via the reproducing property. Such re-interpretation allows generalization of these performance guarantees to COKR. For instance, for a set of $j$ trajectories, let $G_j \in \mathbb{R}^{j\times j}$ denote the Gram matrix of the control occupation kernels, $D_j\in \mathbb{R}^{n\times j}$ denote the end point difference matrix, and $\Gamma_{j}\in\mathbb{R}^{j-1}$ denote the column vector $\begin{pmatrix} \langle \Gamma^{(s)}_{\gamma_{u_1}, u_1}, \Gamma^{(s)}_{\gamma_{u_j}, u_j} \rangle_{H} & \ldots & \langle \Gamma^{(s)}_{\gamma_{u_{j-1}}, u_{j-1}}, \Gamma^{(s)}_{\gamma_{u_j}, u_j} \rangle_{H} \end{pmatrix}^{\top}$. The first $j-1$ trajectories can then be used to predict the $i$-th component of the end point difference for the $j-$th trajectory, denoted by $(D_j)_i^j$, using COKR. The predicted end point difference is given by $(\delta_{j})_i \coloneqq\langle \hat{h}_{i},\Gamma_{\gamma_{u_j},u_j}^{s}\rangle_{H} = (D_{j-1})_i(G_{j-1} + \lambda I_{j-1})^{-1}\Gamma_j$. A straightforward adaptation of the proof of \cite[Theorem 3]{zhdanov2009competing} then establishes the following proposition
\begin{proposition} \label{prop:prediction loss}
    The cumulative end point difference prediction loss satisfies
    \begingroup \small
\begin{multline*}
    \sum_{j=1}^M \frac{((\delta_{j})_i - (D_{j})_i^j)^2}{1+\nu_j} = \\
    \min_{(h)_i \in H}  \left(\sum_{j=1}^M \left[ \left\langle (h)_i, \Gamma^{(s)}_{\gamma_{u_j}, u_j} \right\rangle_H - (D_{j})_i^j\right]^2 + \lambda \| (h)_i \|^2_H\right),
\end{multline*}
\endgroup
where $\nu_j = \frac{\langle \Gamma^{(s)}_{\gamma_{u_j}, u_j}, \Gamma^{(s)}_{\gamma_{u_j}, u_j} \rangle_{H} - \Gamma_j^\top (G_{j-1} + \lambda I_{j-1})^{-1} \Gamma_j}{\lambda} $.
\end{proposition}
\begin{proof}
    Follows from Proposition \ref{prop: FRR prediction loss} below. 
\end{proof}

Furthermore, similar to KRR, it can be shown that in the limit as $M$ goes to infinity, the CORK predictor is at least as good as any continuous function for prediction of end point differences.
\begin{proposition} \label{prop:prediction preformance}
    If $H$ is a vvRKHS of a universal kernel and $X$ is compact then for any continuous function $h:X\to \mathbb{R}^{1\times m+1}$,
    \begin{multline*}
        \lim \sup_{M\to\infty} \frac{1}{M}\Bigg\{\sum_{j=1}^M \left((\delta_{j})_i - (D_{j})_i^j\right)^2 \\
        - \sum_{j=1}^M \left[ \left\langle h, \Gamma^{(s)}_{\gamma_{u_j}, u_j} \right\rangle - (D_{j})_i^j\right]^2\Bigg\} \leq 0
    \end{multline*} 
\end{proposition}
\begin{proof}
    Follows from Proposition \ref{prop: FRR prediction preformance} below.
\end{proof}
\subsection{A generalization of kernel ridge regression}
 It is well-known that the system identification problem can be generalized to other applications where the objective is to approximate an unknown function $g$ using the data $\{(F_j, y_j)\}_{j=1}^{M}$ where $\{F_j\}_{j=1}^{M}\subset H^{*}$, and $H^{*}$ is the set of bounded linear functionals on $H$, such that $[F_j](g) = y_j$ (see, for example, \cite[Theorem 1.3.1]{Wahba1990}). In the functional ridge regression (FRR) approach, the aim is to find a function $\hat{g}\in H$ that solves 
\begin{equation}
    \min_{g\in H} \sum_{i=1}^{M}\left(\left[F_i\right]\left(g\right)-y_i\right)^2+\lambda\|g\|_H^2.\label{eq:FRR_problem}
\end{equation}
The following proposition is a special case of \cite[Theorem 1.3.1]{Wahba1990}. For completeness, a proof is included in the appendix.
\begin{restatable}[]{proposition}{WabhaProposition}
    If $f_i$ is the representative of $F_i$ for all $i$ then there exist constants $ \{c_i\}_{i=1}^M\subseteq \mathbb{R}$ such that the minimizer of \eqref{eq:FRR_problem} is given by $\hat{g} =\sum_{i=1}^M c_i f_i$.
\end{restatable}
\begin{proof}
    See the appendix.
\end{proof}
The vector of constants $ c = \begin{pmatrix} c_1,&\ldots,&c_M \end{pmatrix}^{\top} $ in the FRR minimizer are then given by the resolution of the linear system
\(
    (G_f + \lambda I_M)c = y,
\)
where $y = \begin{pmatrix} y_1,&\ldots,&y_M \end{pmatrix}^{\top}$ denotes the vector of measurements and $G_f\in\mathbb{R}^{M\times M}$ is the Gram matrix of the representing functions $\{f_i\}_{i=1}^M$, whose $i,j-$th element is $\left\langle f_i , f_j \right\rangle_H $. Given a new $F\in H^{*}$, with representative $f\in H$, the action of $F$ on the unknown function $g$ can then be approximated as $[F](g) \approx [F](\hat{g}) = c^{\top} \begin{pmatrix}\left\langle f,f_1 \right\rangle_H,&\ldots,&\left\langle f,f_M \right\rangle_H\end{pmatrix}^{\top}$. Note that for implementation of FRR, all that is required is the ability to compute inner products between the representing functions. 

The COKR problem studied in this paper is an example of FRR where the functionals $F_j$ are the integration functionals $\begin{bmatrix}
    f & g
\end{bmatrix} \mapsto \left\langle\begin{bmatrix}
    f & g
\end{bmatrix}, \Gamma_{\gamma_{u_j},u_j}\right\rangle_H$ and the representatives are the control occupation kernels $\Gamma_{\gamma_{u_j},u_j}$. The KRR problem is another example of FRR, where the functionals $F_j$ are the evaluation functionals $f \mapsto [E_{x_j}](f) = f(x_j)$ and the representatives are data-centered kernels $K(\cdot,x_j)$. In this paper, we show that cumulative error guarantees similar to those derived for KRR in results such as \cite[Theorem 11.1]{vovk2013kernel} also hold for FRR.

\subsection{Cumulative error bounds for FRR}
Let $H$ be a RKHS of $\mathbb{R}^n-$valued continuous functions defined on a compact set $X$. Let $\{F_j\}_{t=1}^\infty$ be a family of bounded linear functionals defined on $\mathcal{C}^0(X,\mathbb{R}^n)$ that satisfies the following assumption.
\begin{assumption}\label{as:1}
There exist positive constants $\overline{F} <\infty$ such that for all $j\in \{1,\ldots,\infty$\}, $f\in \mathcal{C}^0(X,\mathbb{R}^n)$, and $g\in H$, $\left\vert [F_j](f)\right\vert \leq \overline{F} \left\Vert f\right\Vert_\infty$ and $\left\vert [F_j](g)\right\vert \leq \overline{F} \left\Vert g\right\Vert_H$. Furthermore, there exists a positive constant  $\overline{Y} <\infty$ such that for all $j\in \{1,\ldots,\infty$\}, $\left\vert y_j \right\vert \leq \overline{Y}$.
\end{assumption}
If the measurements $y_j$ correspond to the action of the functionals on a continuous scalar field $g$ defined on a compact set, then the existence of $\overline{Y}$ follows from uniform boundedness of the functionals. To facilitate the convergence analysis, let $\hat{g}_j = \arg\min_{g\in H} \sum_{i=1}^{j-1}\left(\left[F_i\right]\left(g\right)-y_i\right)^2+\lambda\|g\|_H^2$ be the FRR minimizer and let $[F_j] 
\left(\hat{g}_j\right)$ be the FRR prediction of the action of the $j-$th functional on the unknown function $g$, both computed using the first $j-1$ functionals and the corresponding $j-1$ measurements.
Note that $[F_j] 
\left(\hat{g}_j\right)$ is a prediction since the measurement $y_j$ is not used to compute $[F_j] 
\left(\hat{g}_j\right)$. Under assumption \ref{as:1}, Proposition \ref{prop:prediction loss} can be generalized to obtain the following result.
\begin{proposition} \label{prop: FRR prediction loss}
    The cumulative prediction loss satisfies
   \begingroup \small $\sum_{j=1}^M \frac{\left(y_j-[F_j] 
\left(\hat{g}_j\right)\right)^2}{1+\nu_j}=
\min_{g\in H} \sum_{j=1}^{M}\left(\left[F_j\right]\left(g\right)-y_j\right)^2+\lambda\|g\|_H^2,
$\endgroup
 where $\nu_j =\frac{\langle f_j,f_j\rangle_H-\overline{f}_j^{\top}\left(G_{f,j-1}+\lambda I\right)^{-1}\overline{f}_j}{\lambda}$, $\overline{f}_j^{\top}=\begin{pmatrix}
    \langle f_1,f_j \rangle_H & \cdots & \langle f_{t-1},f_j \rangle_H
\end{pmatrix}$, and $G_{f,j-1}$ is the Gram matrix of the first $j-1$ representing functions.
\end{proposition}
\begin{proof}
   The proof follows from a straightforward adaptation of the proof of \cite[Theorem 3]{zhdanov2009competing}. 
\end{proof}

Additionally, if $H$ is a vvRKHS of a universal kernel, then the FRR predictor is asymptotically at least as good as any continuous function in the cumulative prediction error metric.
\begin{restatable}[]{proposition}{FRRpredictionpreformance}
 \label{prop: FRR prediction preformance}
    If $H$ is a vvRKHS of a universal kernel, $X$ is compact, and assumption \ref{as:1} is satisfied then for any continuous function $h:X\to \mathbb{R}^{1\times m+1}$
    \begin{multline*}
       \limsup_{M\to \infty}\frac{1}{M}\sum_{j=1}^M \left(\left([F_j] 
\left(\hat{g}_j\right)-y_j\right)^2\!\! -\!\!\left(\left[F_j\right]\left(h\right)-y_j\right)^2\right) \leq 0. 
    \end{multline*} 
\end{restatable}

\begin{proof}
See the appendix.
\end{proof}

Proposition \ref{prop: FRR prediction preformance} shows that the FRR solution converges to a function that preforms better than any continuous function, with respect to the cumulative average prediction error. However, Proposition \ref{prop: FRR prediction preformance} does not imply that the functions $\hat{g}_j$ converge to $g$ in the supremum norm, a limitation FRR shares with KRR. The following section nevertheless demonstrates COKR as an effective tool for approximation of nonlinear dynamical systems.

\section{Numerical Examples}
\subsection{Higher order system - an academic example}
This example utilizes the second order one dimensional nonlinear model of the Duffing oscillator given by 
\(
\ddot{x}=(x-x^3)+(2+\sin(x)) u,
\)
where $f(x)=(x-x^3)$ is the drift function, $g(x)=2+\sin(x)$ is the control effectiveness function, and $u$ is the controller. To approximate the system dynamics, 200 trajectories of the system are recorded, along with the corresponding control signals, starting from a grid of initial conditions $[-3,3] \times [-3,3]$, under a control signal that is composed of the sum of three sinusoidal signals with randomly generated frequencies ranging from $1\ 
 rad/s$ and $3\ rad/s$, with coefficients sampled randomly from $[-1,1]$ . Each trajectory is corrupted by Gaussian measurement noise where the signal-to-noise ratio is about $3\ dB$. The initial velocities are obtained by numerically differentiating the measured noisy trajectories. The recorded trajectories and control signals are then utilized to approximate $f$ and $g$. The kernel used in this example is $K(x,y) = \exp\left(\frac{x^\top y}{\mu}\right)$ with $\mu=5$. 

\textbf{Experiment 1:} The first experiment is an ablation study done to probe the effect of the regularization parameter $\lambda$. To examine the effect of $\lambda$, 22 COKR models of the system are generated from the same dataset of trajectories, each using a different value of $\lambda$. The values of $\lambda$ are selected to be between $10^{-13}$ and $10^7$. Figure \ref{fig:Duffing_COKR_f_tilde_changing_lambda_noise_in_signal} illustrates the results of the ablation study. \begingroup\color{red} The ablation study indicates that there is a wide range of values of $\lambda$ that result in similar model performance. While guidelines to select $\lambda$ directly based on the signal-to-noise ratio are difficult to provide, standard techniques such as cross-validation can be used to numerically identify a value of $\lambda$ that yields acceptable model performance.\endgroup

Based on this ablation study, the regularization parameter is selected to be $\lambda=0.001$ in the following Monte-Carlo trials.
\begin{figure}[h!]
   \centering
    \input{figures/Duffing_COKR_f_tilde_changing_lambda_noise_in_signal}
    \caption{The blue circle marks are the mean of $\left\vert \tilde{f}(x) \right\vert $ and the red square marks are the mean of $\left\vert \tilde{g}(x) \right\vert$ over $x\in[-3,3]$ for different values of $\lambda$ using COKR trained with noisy trajectories.} \label{fig:Duffing_COKR_f_tilde_changing_lambda_noise_in_signal}
\end{figure}

\begin{figure}
   \centering
    \input{figures/Duffing_COKR_f_comparison}
\input{figures/Duffing_COKR_g_comparison}
    \caption{The evaluation of $f$ and $\hat{f}$ (left) and $g$ and $\hat{g}$ (right) for the Duffing oscillator with measurement noise. The solid blue lines represent the true values of the functions and the dotted red lines represent the COKR estimates.} \label{fig:duffing_g_approx_vs_actual}
\end{figure}

\begin{figure}
   \centering
    \input{figures/Duffing_COKR_f_tilda}
    \input{figures/Duffing_COKR_g_tilda}
    \caption{The evaluation of $\left\vert \tilde{f} \right\vert$ (left) and $\left\vert \tilde{g} \right\vert$ (right) for the Duffing oscillator with measurement noise.} \label{fig:g_tilde}
\end{figure}



\textbf{Experiment 2:} A Monte-Carlo simulation with 1000 trials is conducted to evaluate robustness of the developed method to excitation signals and sensor noise. Figure \ref{fig:duffing_g_approx_vs_actual} shows the estimated and actual values of the drift function and 
the estimated and actual values of the control effectiveness function in the first Monte-Carlo trial. Figure \ref{fig:g_tilde} shows the drift approximation error, $\tilde{f}$, as a function of $x$ and 
the control effectiveness approximation error, $\tilde{g}$, as a function of $x$ in the first Monte-Carlo trial. For each trial, the maximum values of $\left\vert \tilde{f} \right\vert$ over the interval $[-3,3]$ are collected using COKR and SINDYc for comparison. Box plots for the Monte-Carlo trials are presented in Figure 
\ref{fig:g_error_boxPolt} for the maximum errors  $\max_{x\in[-3,3]}\left\vert \tilde{f} \right\vert$ and $\max_{x\in[-3,3]}\left\vert \tilde{g} \right\vert$, respectively. 



\begin{figure}
   \centering
    \pgfplotstableread[col sep=comma,header=false] {data/error_by_iteration_f_both_methods.dat}\Table
    \input{drawBoxPlots}
    \hspace{1em}
    \pgfplotstableread[col sep=comma,header=false] {data/error_by_iteration_g_both_methods.dat}\Table
    \input{drawBoxPlots}
    \caption{Monte-Carlo results of the max of $\left\vert \tilde{f} \right\vert$ (left) and $\left\vert \tilde{g} \right\vert$ (right) over 1000 trials.} \label{fig:g_error_boxPolt}
\end{figure}

\subsection{Control-affine model of a two-link robot manipulator}

This example utilizes a first-order four-dimensional model of a two-link robot manipulator given by
\begingroup \small
\begin{equation}\label{eq:Full_State_Two_link_dynamics}
\begin{pmatrix} \dot{x}_1 \\ \dot{x}_2\end{pmatrix}=\begin{pmatrix} x_2 \\ -M^{-1}(x_1) C(x_1,x_2)\end{pmatrix} +
\begin{pmatrix} 0_{2\times 2} \\ M^{-1}(x_1) \end{pmatrix} \begin{pmatrix} \tau_1 \\ \tau_2 \end{pmatrix},
\end{equation}
\endgroup
where $x_1:=\begin{pmatrix}q_1\\q_2\end{pmatrix}$, $x_2:=\begin{pmatrix}\dot{q}_1\\\dot{q}_2\end{pmatrix}$, $q_1$, and $q_2$ are the angular positions of the two links, respectively, $M(x_1)=\begin{pmatrix}  p_1+2 p_3 \cos(q_2) & p_2+p_3 \cos(q_2) \\ p_2+p_3 \cos(q_2) & p_2 \end{pmatrix}$ is the inertia matrix, $C(x_1,x_2)=\left(V(x_1,x_2) + F_d \right) x_2 + F_s(x_2)$, where $F_d=\begin{pmatrix}  f_{d_1} & 0 \\ 0 & f_{d_2} \end{pmatrix}$ denotes the viscous friction at the two joints, $F_s(x_2)=\begin{pmatrix}  f_{s_1} \tanh(\dot{q}_1)  \\ f_{s_2} \tanh(\dot{q}_2) \end{pmatrix}$ denotes the static friction at the two joints,  $V (x_1,x_2) = \begin{pmatrix} -p_3  \sin(q_2)  \dot{q}_2 & -p_3  \sin(q_2) (\dot{q}_1+\dot{q}_2) \\ p_3  \sin(q_2) \cdot \dot{q}_1 & 0 \end{pmatrix}$ is the centrifugal-Coriolis matrix and $\tau_1$ and $\tau_2$ are the control torques applied by the two motors. The parameters are given by $p_1=3.473$, $p_2=0.196$,  $p_3=0.242$, $f_{d_1}=5.3$, $f_{d_2}=1.1$, $f_{s_1}=8.45$, and $f_{s_2}=2.35$.


\textbf{Experiment 1: } Trajectories of the system are collected, starting from 1000 different initial conditions, sampled using pseudorandom Halton sampling over a cube of side 6 centered at the origin of $\mathbb{R}^4$. Each trajectory is recorded using control signals $\tau_1$ and $\tau_2$, comprised of a sum of three sinusoidal signals with randomly generated frequencies and coefficients. For each run, the control signals and the trajectories of the system are recorded to generate the database that is used to approximate $f$ and $g$. Since this system has two controllers $\tau_1$ and $\tau_2$, the control effectiveness matrix $g$ is decomposed into $g_1$ and $g_2$, where $g_1$ is the first row of $g$ and $g_2$ is the second row of $g$.

\begin{figure}
   \centering
    \input{figures/TwoLink_f_tilda}

    \input{figures/TwoLink_g1_tilda}
    \input{figures/TwoLink_g2_tilda}
    \caption{The relative errors $ \frac{\left\Vert \tilde{f}(x) \right\Vert}{\left\Vert f(x) \right\Vert}$ (top), $\frac{\Vert \tilde{g}_1(x) \Vert}{\Vert g_1(x) \Vert}$ (bottom left), and $\frac{\Vert \tilde{g}_2(x) \Vert}{\Vert g_2(x) \Vert}$ (bottom right) evaluated at 100 points, indexed by decreasing distance from the origin.} \label{fig:Two_link_g2_tilde}
\end{figure}

\begin{figure}
   \centering
    \input{figures/TwoLink_torque_true}
    \input{figures/TwoLink_torque_estimated}
    \caption{Regulation result using the true model (left) and estimated model (right) for computed torque control. The green line is the the angular position $q_1 (rad)$, the dashed red line is the the angular position $q_2 (rad)$, the dashed black line is the the angular velocity $\dot{q}_1 (rad/s)$, and the dashed blue line is the the angular velocity $\dot{q}_2 (rad/s)$.} \label{fig:computed_torque_regulator}
\end{figure}


The approximation error between the identified system and the actual system is computed through evaluation at 100 sample points sampled using pseudorandom Halton sampling over a cube of side 4 centered at the origin of $\mathbb{R}^4$. The functions $f$, $\hat{f}$, $g$, and $\hat{g}$ are evaluated at each vertex of the grid to yield the approximation errors $\tilde{f}$, $\tilde{g_1}$, and $\tilde{g_2}$.

Figure \ref{fig:Two_link_g2_tilde} shows the relative errors $ \frac{\left\Vert \tilde{f}(x) \right\Vert}{\left\Vert f(x) \right\Vert}$, $ \frac{\left\Vert \tilde{g}_1(x) \right\Vert}{\left\Vert g_1(x) \right\Vert}$, and $ \frac{\left\Vert \tilde{g}_2(x) \right\Vert}{\left\Vert g_2(x) \right\Vert}$, respectively as a function of the node index, \begingroup \color{red} all using a regularization parameter of $\lambda=0.001$.  \endgroup 
 
\textbf{Experiment 2: }To demonstrate the utility of the model developed via the novel system identification method, we simulate the actual two-link robot manipulator using a fourth order Runge-Kutta method with a computed torque controller designed to regulate the angular positions and velocities to zero. The computed torque controller is calculated using the estimated model, and is given by 
\begin{equation}\label{eq:computed_toruqe_regulator}
\tau = \hat{M}(q) [-K_v \dot{q}- K_p q]+\hat{C}(q,\dot{q})
\end{equation}
where $\hat{M}(q)$ and $\hat{C}(q,\dot{q})$ are recovered from the estimated data-driven drift and control-effectiveness functions $\hat{f}$ and $\hat{g}$, and $K_p=\begin{pmatrix}20 & 0 \\ 0 & 20 \end{pmatrix}$, 
$K_v=\begin{pmatrix}30 & 0 \\ 0 & 30 \end{pmatrix}$ are feedback gains.
Figure \ref{fig:computed_torque_regulator} shows the performance of the computed torque regulator used on the true two-link robot manipulator model.

\section{Discussion}

In the first numerical example, where the Duffing oscillator system is identified, the approximation  $\hat{f}$ is nearly identical to the true ${f}$ as seen in Figure \ref{fig:duffing_g_approx_vs_actual} and the maximum error is $0.61$, which is a $2.6\%$ error as seen in Figure \ref{fig:g_tilde}. The approximation $\hat{g}$ captures the underlining structure of $g$ where $g(x)=2+\sin(x)$, but $\hat{g}$ deviates slightly from the true value as seen in Figure \ref{fig:duffing_g_approx_vs_actual} with a maximum error of $0.4$  which is a $23\%$ error as seen in Figure \ref{fig:g_tilde}. The maximum errors are a function of kernel type, parameters, number and spatial coverage of the recorded data, and measurement noise, and can be reduced through segmentation of the recorded trajectories to generate more data points and increase the resolution of the approximations.

\begin{table}
\resizebox{\columnwidth}{!}{%
\begin{tabular}{ccc|cc}
\multicolumn{1}{l}{}                                                             & \multicolumn{2}{c|}{COKR}            & \multicolumn{2}{c}{SIDNYc}            \\ \hline
\multicolumn{1}{l|}{}                                                            & \multicolumn{1}{c|}{Mean}   & SD     & \multicolumn{1}{c|}{Mean}   & SD      \\
\multicolumn{1}{c|}{max$  \left(  \left\vert \tilde{f}(x) \right\vert \right) $} & \multicolumn{1}{c|}{1.0322} & 0.1765 & \multicolumn{1}{c|}{2.5224} & 13.3437 \\
\multicolumn{1}{c|}{}                                                            & \multicolumn{1}{c|}{}       &        & \multicolumn{1}{c|}{}       &         \\
\multicolumn{1}{c|}{max$  \left(  \left\vert \tilde{g}(x) \right\vert \right) $} & \multicolumn{1}{c|}{0.4718} & 0.2346 & \multicolumn{1}{c|}{1.1095} & 0.3417 
\end{tabular}%
}
\caption{Comparison result of the Monte-Carlo simulations of the two methods.} \label{tab:comparison}
\end{table}

The Monte-Carlo trials summarized in Figure 
\ref{fig:g_error_boxPolt} demonstrate the robustness of the developed system identification method to  excitation signals and measurement noise. Figure 
\ref{fig:duffing_g_approx_vs_actual} show the full error plots of a representative sample trial.

Figure \ref{fig:g_error_boxPolt} indicates that COKR is less sensitive to measurement noise than SIDNYc. In Table \ref{tab:comparison} the mean and standard deviation of max$  \left(  \left\vert \tilde{f}(x) \right\vert \right) $ and max$  \left(  \left\vert \tilde{g}(x) \right\vert \right) $ for COKR and SIDNYc are listed, which shows that COKR is more precise than SIDNYc. This result is expected in a scenario where measurement noise is present, since COKR uses numerical integration whereas SIDNYc uses numerical differentiation. Although the median of the error measurement of SINDYc for $f$ in Figure \ref{fig:g_error_boxPolt} is lower than that of COKR, SINDYc results in a significant number of outliers, where the worst performing run produces an error measurement that is two orders of magnitude higher than COKR.

In the second numerical example, where the model of a two-link robot manipulator is identified, Figure 
\ref{fig:Two_link_g2_tilde} indicates that  the estimation errors get worse as one approaches the boundary of the domain covered by the trajectories, and better as one approaches the origin. While the errors are hard to visualize as a function of distance from the training data in four dimensions, the trends observed in Figure 
\ref{fig:Two_link_g2_tilde} could heuristically be attributed to limited coverage of the corners of the domain. 

The identified system $\begin{pmatrix} \hat{f} & \hat{g} \end{pmatrix}$ is used to compute the necessary torque to regulate the actual two-link robot manipulator, which results in a performance similar to a computed torque controller implemented using exact model knowledge , as seen in Figure \ref{fig:computed_torque_regulator}. 
The main downside of regulating the system using $\begin{pmatrix} \hat{f} & \hat{g} \end{pmatrix}$ is the long computation time necessary to evaluate $\begin{pmatrix} \hat{f} & \hat{g} \end{pmatrix}$ at any given $x$, which highlights the need for optimization of the evaluation function to improve computational performance.

\section{Conclusion}\label{sec:conclusion}

A data-driven control occupation kernel regression method is developed in this manuscript for identification of nonlinear affine control systems, where identification of higher order systems is possible without numerical differentiation for higher order systems. As a result, as indicated by the numerical experiments, the developed method is robust to measurement noise. While the model developed using this method can be used to compute the feedforward component of a controller, depending on the number of trajectories used for modeling, evaluation of the model at a given state can be computationally expensive. Further research is required to develop a more efficient evaluation method to render it useful for real-time feedback. 

A distinct advantage of the use of an occupation kernel basis to approximate the dynamics of the system is that under the assumptions of the framework developed in this paper, the solution of the minimization problem in \eqref{eq:regression} is guaranteed to be a linear combination of the occupation kernels by the representer theorem. This means that the occupation kernels are natural basis functions arising from a dynamics context, and in principle they should perform better than less structured bases. In particular, through the representer theorem, it is clear that the occupation kernel basis will, in general, result in models that yield a lower modeling error in the Hilbert space norm than any other generic basis, resulting in better generalization of the estimates away from the training data. \begingroup\color{red} Data-richness in the developed algorithm is related to orthogonality of the control occupation kernels in the Hilbert space. Examination of the relationship between excitation in the input signal in terms of independent frequencies, and orthogonality of the resulting control occupation kernels, is a part of future research.\endgroup

\appendix

\densityProposition*
\begin{proof}
    Select $\gamma_0$ such that $h(\gamma_0) \neq 0$. Continuity of $h$ and $\gamma_{u,\gamma_0}$ can then be invoked to conclude that for any constant control signal $u(t) = b$ such that $h(\gamma_0) \begin{pmatrix} 1 \\ b \end{pmatrix} \neq 0$, there exists a $T_u > 0$ for which $ \int_{0}^{T_u} h(\gamma_{u,\gamma_0}(t)) \begin{pmatrix} 1 \\ b \end{pmatrix} \mathrm{d}t \neq 0$. For example, one can select $T_u = \min\left\{T,\inf_t\left\{h(\gamma_{u,\gamma_0}(t))\begin{pmatrix} 1 \\ b \end{pmatrix}=0\right\}\right\}$. A straightforward extension of the above argument to the higher order case allows us to conclude that for any order $s$ and any nonzero $h$, there exist $\gamma_0$, $u$, and $T_u$ such that $\left\langle h , \Gamma_{\gamma_{u,\gamma_0},u}^{(s)} \right\rangle_H \neq 0$. That is, $H\cap A_s^\perp = \{0\}$, and as a result, $(A_s^\perp)^\perp = H$. Since $(A_s^\perp)^\perp = \overline{\vspan A_s}$, we conclude that $H = \overline{\vspan A_s}$.
\end{proof}
\WabhaProposition*
\begin{proof}
write $g= g_1+g_2$, with $g_1\in \vspan\{f_i\}_{i=1}^M \coloneqq S$ and $g_2 \in S^{\perp}$. Then, the first term in the cost,  $\sum_{i=1}^{M}\left(\left[F_i\right]\left(g\right)-y_i\right)^2= \sum_{i=1}^{M}\left(\langle g_1,f_i\rangle+\langle g_2,f_i\rangle -y_i\right)^2=\sum_{i=1}^{M}\left(\langle g_1,f_i\rangle -y_i\right)^2$, is independent of $g_2$. The second term is given by 
$\|g\|_H^2=\langle g_1+g_2,g_1+g_2\rangle_H=\langle g_1,g_1\rangle_H+2\langle g_1,g_2 \rangle_H+\langle g_2,g_2\rangle=\|g_1\|_H^2+\|g_2\|_H^2$.
Since the first term is independent of $g_2$ and the second term is monotonic in $\|g_2\|_H^2$, for any $\hat{g}$ that minimizes the cost, $\|g_2\|_H^2=0$.
\end{proof}
\FRRpredictionpreformance*

\begin{proof}
Given $\epsilon > 0$ and a continuous function $g$, let $\overline{g}_{\epsilon}\in H$ be a function that satisfies $\sup_{x\in X} \left\Vert g(x) - \overline{g}_{\epsilon}(x)\right\Vert\leq \epsilon$. The probabilistic interpretation in section 5 of \cite{zhdanov2009competing} implies that $\langle f_j,f_j\rangle_H-\overline{f}_j^{\top}\left(G_{f,j-1}+\lambda I\right)^{-1}\overline{f}_j > 0$ and as a result, Proposition \ref{prop: FRR prediction loss} implies that
$\sum_{j=1}^M \left(y_j-[F_j] 
\left(\hat{g}_j\right)\right)^2 < \min_{h\in H} \sum_{j=1}^{M}\left(\left[F_j\right]\left(h\right)-y_j\right)^2+\lambda\|h\|_H^2.$ In particular, with $h = \overline{g}_{\epsilon}$, 
    $\sum_{j=1}^M \left(y_j-[F_j] 
\left(\hat{g}_j\right)\right)^2 - \sum_{j=1}^{M}\left(\left[F_j\right]\left(g\right)-y_j\right)^2 
+ \sum_{j=1}^{M}\left(\left(\left[F_j\right]\left(g\right)-y_j\right)^2 -\left(\left[F_j\right]\left(\overline{g}_{\epsilon}\right)-y_j\right)^2 \right)< \lambda\|\overline{g}_{\epsilon}\|_H^2$,
Expanding the squares in the last term on the left hand side, using linearity of $ F_j$  and the identity $ \left(\left[F_j\right]\left(g\right)\right)^2-\left(\left[F_j\right]\left(\overline{g}_{\epsilon}\right)\right)^2=\left(\left[F_j\right]\left(g+\overline{g}_{\epsilon}\right)\right)\left(\left[F_j\right]\left(g-\overline{g}_{\epsilon}\right)\right)$, and using 
assumption \ref{as:1},
$
    \sum_{j=1}^M \left(y_j-[F_j] \left(\hat{g}_j\right)\right)^2 - \sum_{j=1}^{M}\left(\left[F_j\right]\left(g\right)-y_j\right)^2  < \lambda\|\overline{g}_{\epsilon}\|_H^2 + \sum_{j=1}^{M}\overline{F}^2\left\Vert g+\overline{g}_{\epsilon}\right\Vert_\infty \left\Vert g-\overline{g}_{\epsilon}\right\Vert_\infty + 2\overline{F}\left\Vert \overline{g}_{\epsilon}-g \right\Vert_\infty \overline{Y}
$
Adding and subtracting $g$ to the first term in the second row,
\begin{multline*}  
\sum_{j=1}^M \left(y_j-[F_j] 
\left(\hat{g}_j\right)\right)^2 - \sum_{j=1}^{M}\left(\left[F_j\right]\left(g\right)-y_j\right)^2  < \lambda\|\overline{g}_{\epsilon}\|_H^2 +\\ \sum_{j=1}^{M}\left(\overline{F}^2 2\left\Vert g\right\Vert_\infty + 2\overline{F}\overline{Y}\right)\left\Vert \overline{g}_{\epsilon}-g \right\Vert_\infty + \left\Vert \overline{g}_{\epsilon}-g\right\Vert_\infty^2.
\end{multline*}
Since $\sup_{x\in X} \left\Vert g(x) - \overline{g}_{\epsilon}(x)\right\Vert\leq \epsilon$,

$\frac{1}{M}\left(\sum_{j=1}^M \left(y_j-[F_j] 
\left(\hat{g}_j\right)\right)^2 - \sum_{j=1}^{M}\left(\left[F_j\right]\left(g\right)-y_j\right)^2\right)  < \\\frac{\lambda}{M}\|\overline{g}_{\epsilon}\|_H^2+ \left(\overline{F}^2\left\Vert 2g\right\Vert_\infty+ 2\overline{F}\overline{Y}\right){\epsilon} + {\epsilon}^2.$ 


Taking the limit superior as $M \to \infty$, and then taking the limit as $\epsilon \to 0$, the proposition is established.
\end{proof}




\small
\bibliographystyle{IEEEtran}
\bibliography{references}

\end{document}

%% file: figures/Duffing_COKR_f_tilde_changing_lambda_noise_in_signal.tex
\begin{tikzpicture}
    \begin{loglogaxis}[
        xlabel={ $\lambda$ },
        ylabel={  },
        legend pos = outer north east,
        legend style={nodes={scale=0.5, transform shape}},
        enlarge y limits=0.1,
        enlarge x limits=0,
        height = 0.4\columnwidth,
        width = \columnwidth,
        label style={font=\scriptsize},
        tick label style={font=\scriptsize}
    ]
        \addplot [only marks, blue] table [x index=0, y index=1]{data/Duffing_COKR_f_g_tilde_changing_lambda_noise_in_signal.dat};

        \addplot [only marks, red, mark=square*] table [x index=0, y index=2]{data/Duffing_COKR_f_g_tilde_changing_lambda_noise_in_signal.dat};
        
    \end{loglogaxis}
\end{tikzpicture}

%% file: figures/Duffing_COKR_f_comparison.tex
\begin{tikzpicture}
    \begin{axis}[
        xlabel={$x$},
        ylabel={$f(x)$},
        enlarge y limits=0.1,
        enlarge x limits=0,
        height = 0.3\columnwidth,
        width = 0.5\columnwidth,
        label style={font=\scriptsize},
        tick label style={font=\scriptsize},
        ylabel shift = -8 pt,
    ]
        \addplot [thick, blue] table [x index=0, y index=1]{data/Duffing_COKR_f.dat};
        \addplot [dashed, red] table [x index=0, y index=2]{data/Duffing_COKR_f.dat};
   
    \end{axis}
\end{tikzpicture}

%% file: figures/Duffing_COKR_g_comparison.tex
\begin{tikzpicture}
    \begin{axis}[
        xlabel={ $x$ },
        ylabel={  $g(x)$},
        enlarge y limits=0.1,
        enlarge x limits=0,
        height = 0.3\columnwidth,
        width = 0.5\columnwidth,
        label style={font=\scriptsize},
        tick label style={font=\scriptsize}
    ]
        \addplot [thick, blue] table [x index=0, y index=1]{data/Duffing_COKR_g.dat};
        \addplot [dashed, red] table [x index=0, y index=2]{data/Duffing_COKR_g.dat};
   
    \end{axis}
\end{tikzpicture}

%% file: figures/Duffing_COKR_f_tilda.tex
\begin{tikzpicture}
    \begin{axis}[
        xlabel={ $x$ },
        ylabel={  $\tilde{f}(x)$},
        legend pos = outer north east,
        legend style={nodes={scale=0.5, transform shape}},
        enlarge y limits=0.1,
        enlarge x limits=0,
        height = 0.3\columnwidth,
        width = 0.5\columnwidth,
        label style={font=\scriptsize},
        tick label style={font=\scriptsize}
    ]
        \addplot [thick, blue] table [x index=0, y index=3]{data/Duffing_COKR_f.dat};

    \end{axis}
\end{tikzpicture}

%% file: figures/Duffing_COKR_g_tilda.tex
\begin{tikzpicture}
    \begin{axis}[
        xlabel={ $x$},
        ylabel={ $\tilde{g}(x)$},
        legend pos = outer north east,
        legend style={nodes={scale=0.5, transform shape}},
        enlarge y limits=0.1,
        enlarge x limits=0,
        height = 0.3\columnwidth,
        width = 0.5\columnwidth,
        label style={font=\scriptsize},
        tick label style={font=\scriptsize}
    ]
        \addplot [thick, blue] table [x index=0, y index=3]{data/Duffing_COKR_g.dat};

    \end{axis}
\end{tikzpicture}

%% file: drawBoxPlots.tex
\pgfplotsset{
    boxplot/lower notch/.initial=\pgfutil@empty,
    boxplot/upper notch/.initial=\pgfutil@empty,
    boxplot/notch width/.initial=0.5,
    boxplot/draw/box/.code={%
        \draw[/pgfplots/boxplot/every box/.try]
        (boxplot box cs:\pgfplotsboxplotvalue{lower quartile},0)
        -- (boxplot box cs:\pgfplotsboxplotvalue{lower notch},0)
        -- (boxplot box cs:\pgfplotsboxplotvalue{median},0.5-\pgfplotsboxplotvalue{notch width}/2)
        -- (boxplot box cs:\pgfplotsboxplotvalue{upper notch},0)
        -- (boxplot box cs:\pgfplotsboxplotvalue{upper quartile},0)
        -- (boxplot box cs:\pgfplotsboxplotvalue{upper quartile},1)
        -- (boxplot box cs:\pgfplotsboxplotvalue{upper notch},1)
        -- (boxplot box cs:\pgfplotsboxplotvalue{median},0.5+\pgfplotsboxplotvalue{notch width}/2)
        -- (boxplot box cs:\pgfplotsboxplotvalue{lower notch},1)
        -- (boxplot box cs:\pgfplotsboxplotvalue{lower quartile},1)
        -- cycle;
    },%
    boxplot/draw/median/.code={%
        \draw[/pgfplots/boxplot/every median/.try, color=ForestGreen]
        (boxplot box cs:\pgfplotsboxplotvalue{median},0.5-\pgfplotsboxplotvalue{notch width}/2)
        -- (boxplot box cs:\pgfplotsboxplotvalue{median},0.5+\pgfplotsboxplotvalue{notch width}/2)
        ;
    },%
    boxplot prepared from table/.code={
        \def\tikz@plot@handler{\pgfplotsplothandlerboxplotprepared}%
        \pgfplotsset{
            /pgfplots/boxplot prepared from table/.cd,
            #1,
        }
    },
    /pgfplots/boxplot prepared from table/.cd,
    table/.code={\pgfplotstablecopy{#1}\to\boxplot@datatable},
    row/.initial=0,
    make style readable from table/.style={
        #1/.code={
            \pgfplotstablegetelem{\pgfkeysvalueof{/pgfplots/boxplot prepared from table/row}}{##1}\of\boxplot@datatable
            \pgfplotsset{boxplot/#1/.expand once={\pgfplotsretval}}
        }
    },
    make style readable from table=lower whisker,
    make style readable from table=upper whisker,
    make style readable from table=lower quartile,
    make style readable from table=upper quartile,
    make style readable from table=median,
    make style readable from table=lower notch,
    make style readable from table=upper notch
}
\begin{tikzpicture}
    \pgfplotstablegetrowsof{\Table}
    \pgfmathsetmacro\lastRowIndx{\pgfplotsretval-1}
    \pgfplotstablegetcolsof{\Table}
    \pgfmathsetmacro\lastcolIndx{\pgfplotsretval-1}
    \edef\lastCol{\pgfplotsretval}
    \ifnum \lastCol > 8
        \def\outlierList{8}
        \pgfplotsinvokeforeach{9,...,\lastcolIndx}{\edef\outlierList{\outlierList,#1}}
        \pgfplotstabletranspose[columns/.expanded=\outlierList]{\Outliers}{\Table}
    \fi
    \begin{axis}[
        width = 0.45\columnwidth,
        height = 0.5\columnwidth,
        label style={font=\scriptsize},
        tick label style={font=\scriptsize},
        boxplot/draw direction=y,
        ymode = log,
        xtick = {0,...,\lastRowIndx},
        xticklabel style = {align=center, font=\small, rotate=60},
        xticklabels from table={\Table}{[index] 0}
        ]
        \pgfplotsinvokeforeach {0,...,\lastRowIndx} {
            \ifnum \lastCol > 8
                \addplot[
                    mark=+,
                    mark options={red},
                    boxplot prepared from table={
                        table=\Table,
                        row=#1,
                        lower whisker=1,
                        upper whisker=5,
                        lower quartile=2,
                        upper quartile=4,
                        lower notch=6,
                        upper notch=7,
                        median=3,
                    },
                    boxplot prepared={draw position=#1}
                ] table [y=#1] {\Outliers};
            \else
                \addplot[
                    mark=+,
                    mark options={red},
                    boxplot prepared from table={
                        table=\Table,
                        row=#1,
                        lower whisker=1,
                        upper whisker=5,
                        lower quartile=2,
                        upper quartile=4,
                        lower notch=6,
                        upper notch=7,
                        median=3,
                    },
                    boxplot prepared={draw position=#1}
                ] coordinates {};
            \fi
        }
    \end{axis}
\end{tikzpicture}

%% file: figures/TwoLink_f_tilda.tex
\begin{tikzpicture}
    \begin{axis}[
        xlabel={ Index },
        ylabel={  $\frac{\Vert \tilde{f}(x) \Vert}{\Vert f(x) \Vert}$},
        legend pos = outer north east,
        legend style={nodes={scale=0.5, transform shape}},
        enlarge y limits=0.1,
        enlarge x limits=0,
        height = 0.3\columnwidth,
        width = 0.5\columnwidth,
        label style={font=\scriptsize},
        tick label style={font=\scriptsize},
        ylabel shift = -0 pt,
    ]
        \addplot [only marks, blue, mark size=1pt] table [x index=0, y index=1]{data/TwoLink_f_tilda.dat};

    \end{axis}
\end{tikzpicture}

%% file: figures/TwoLink_g1_tilda.tex
\begin{tikzpicture}
    \begin{axis}[
        xlabel={ Index },
        ylabel={  $\frac{\Vert \tilde{g}_1(x) \Vert}{\Vert g_1(x) \Vert}$},
        legend pos = outer north east,
        legend style={nodes={scale=0.5, transform shape}},
        enlarge y limits=0.1,
        enlarge x limits=0,
        height = 0.3\columnwidth,
        width = 0.45\columnwidth,
        label style={font=\scriptsize},
        tick label style={font=\scriptsize},
        ylabel shift = -5 pt,
    ]
        \addplot [only marks, blue, mark size=1pt] table [x index=0, y index=1]{data/TwoLink_g1_tilda.dat};

    \end{axis}
\end{tikzpicture}

%% file: figures/TwoLink_g2_tilda.tex
\begin{tikzpicture}
    \begin{axis}[
        xlabel={ Index },
        ylabel={  $\frac{\Vert \tilde{g}_2(x) \Vert}{\Vert g_2(x) \Vert}$},
        legend pos = outer north east,
        legend style={nodes={scale=0.5, transform shape}},
        enlarge y limits=0.1,
        enlarge x limits=0,
        height = 0.3\columnwidth,
        width = 0.45\columnwidth,
        label style={font=\scriptsize},
        tick label style={font=\scriptsize},
        ylabel shift = -5 pt,
    ]
        \addplot [only marks, blue, mark size=1pt] table [x index=0, y index=1]{data/TwoLink_g2_tilda.dat};

    \end{axis}
\end{tikzpicture}

%% file: figures/TwoLink_torque_true.tex
\begin{tikzpicture}
    \begin{axis}[
        xlabel={ Time (s) },
        legend pos = north east,
        legend style={nodes={scale=0.4, transform shape}},
        enlarge y limits=0.1,
        enlarge x limits=0,
        height = 0.3\columnwidth,
        width = 0.5\columnwidth,
        label style={font=\scriptsize},
        tick label style={font=\scriptsize}
    ]
        \addplot [thick,ForestGreen] table [x index=0, y index=1]{data/TwoLink_torque_true.dat};
        \addplot [thick,dashed,red] table [x index=0, y index=2]{data/TwoLink_torque_true.dat};
        \addplot [thick,dotted,black] table [x index=0, y index=3]{data/TwoLink_torque_true.dat};
        \addplot [thick,dash dot,blue] table [x index=0, y index=4]{data/TwoLink_torque_true.dat};

    \end{axis}
\end{tikzpicture}

%% file: figures/TwoLink_torque_estimated.tex
\begin{tikzpicture}
    \begin{axis}[
        xlabel={ Time (s) },
        legend pos = north east,
        legend style={nodes={scale=0.4, transform shape}},
        enlarge y limits=0.1,
        enlarge x limits=0,
        height = 0.3\columnwidth,
        width = 0.5\columnwidth,
        label style={font=\scriptsize},
        tick label style={font=\scriptsize}
    ]
        \addplot [thick,ForestGreen] table [x index=0, y index=1]{data/TwoLink_torque_estimated.dat};
        \addplot [thick,dashed,red] table [x index=0, y index=2]{data/TwoLink_torque_estimated.dat};
        \addplot [thick,dotted,black] table [x index=0, y index=3]{data/TwoLink_torque_estimated.dat};
        \addplot [thick,dash dot,blue] table [x index=0, y index=4]{data/TwoLink_torque_estimated.dat};

    \end{axis}
\end{tikzpicture}

%% file: root.bbl
\begin{thebibliography}{10}
\providecommand{\url}[1]{#1}
\csname url@samestyle\endcsname
\providecommand{\newblock}{\relax}
\providecommand{\bibinfo}[2]{#2}
\providecommand{\BIBentrySTDinterwordspacing}{\spaceskip=0pt\relax}
\providecommand{\BIBentryALTinterwordstretchfactor}{4}
\providecommand{\BIBentryALTinterwordspacing}{\spaceskip=\fontdimen2\font plus
\BIBentryALTinterwordstretchfactor\fontdimen3\font minus \fontdimen4\font\relax}
\providecommand{\BIBforeignlanguage}[2]{{%
\expandafter\ifx\csname l@#1\endcsname\relax
\typeout{** WARNING: IEEEtran.bst: No hyphenation pattern has been}%
\typeout{** loaded for the language `#1'. Using the pattern for}%
\typeout{** the default language instead.}%
\else
\language=\csname l@#1\endcsname
\fi
#2}}
\providecommand{\BIBdecl}{\relax}
\BIBdecl

\bibitem{SCC.Chowdhary.Johnson2011}
G.~V. Chowdhary and E.~N. Johnson, ``Theory and flight-test validation of a concurrent-learning adaptive controller,'' \emph{J. Guid. Control Dynam.}, vol.~34, no.~2, pp. 592--607, Mar. 2011.

\bibitem{SCC.Brunton.Proctor.ea2016a}
S.~L. Brunton, J.~L. Proctor, and J.~N. Kutz, ``Discovering governing equations from data by sparse identification of nonlinear dynamical systems,'' \emph{Proc. Nat. Acad. Sci. U.S.A.}, vol. 113, no.~15, pp. 3932--3937, 2016.

\bibitem{SCC.Ljung1998}
L.~Ljung, ``System identification,'' in \emph{Signal Analysis and Prediction}, ser. Applied and Numerical Harmonic Analysis, A.~Proch{\'{a}}zka, J.~Uhlir, P.~W.~J. Rayner, and N.~G. Kingsbury, Eds.\hskip 1em plus 0.5em minus 0.4em\relax Birkh{\"{a}}user Boston, 1998, pp. 163--173.

\bibitem{SCC.Brunton.Proctor.ea2016}
S.~L. Brunton, J.~L. Proctor, and J.~N. Kutz, ``Sparse identification of nonlinear dynamics with control {(SINDYc)},'' \emph{IFAC-PapersOnLine}, vol.~49, no.~18, pp. 710--715, 2016.

\bibitem{pillonetto2014kernel}
G.~Pillonetto, F.~Dinuzzo, T.~Chen, G.~De~Nicolao, and L.~Ljung, ``Kernel methods in system identification, machine learning and function estimation: A survey,'' \emph{Automatica}, vol.~50, no.~3, pp. 657--682, 2014.

\bibitem{pillonetto2010new}
G.~Pillonetto and G.~De~Nicolao, ``A new kernel-based approach for linear system identification,'' \emph{Automatica}, vol.~46, no.~1, pp. 81--93, 2010.

\bibitem{pillonetto2011new}
G.~Pillonetto, M.~H. Quang, and A.~Chiuso, ``A new kernel-based approach for nonlinearsystem identification,'' \emph{IEEE Transactions on Automatic Control}, vol.~56, no.~12, pp. 2825--2840, 2011.

\bibitem{pillonetto2015tuning}
G.~Pillonetto and A.~Chiuso, ``Tuning complexity in regularized kernel-based regression and linear system identification: The robustness of the marginal likelihood estimator,'' \emph{Automatica}, vol.~58, pp. 106--117, 2015.

\bibitem{SCC.Rosenfeld.Russo.ea2024}
J.~A. Rosenfeld, B.~Russo, R.~Kamalapurkar, and T.~Johnson, ``The occupation kernel method for nonlinear system identification,'' \emph{SIAM J. Control Optim.}, to appear, see {arXiv:1909.11792}.

\bibitem{SCC.Russo.Kamalapurkar.ea2022}
B.~P. Russo, R.~Kamalapurkar, D.~Chang, and J.~A. Rosenfeld, ``Motion tomography via occupation kernels,'' \emph{J. Comput. Dyn.}, vol.~9, no.~1, pp. 27--45, 2022.

\bibitem{SCC.Li.Rosenfeld2020}
X.~Li and J.~A. Rosenfeld, ``Fractional order system identification with occupation kernel regression,'' \emph{IEEE Control Syst. Lett.}, to appear.

\bibitem{SCC.Rosenfeld.Kamalapurkar.ea2022}
J.~A. Rosenfeld, R.~Kamalapurkar, L.~F. Gruss, and T.~T. Johnson, ``Dynamic mode decomposition for continuous time systems with the {L}iouville operator,'' \emph{J. Nonlinear Sci.}, vol.~32, no.~1, pp. 1--30, Feb. 2022.

\bibitem{SCC.Lasserre.Henrion.ea2008}
J.~B. Lasserre, D.~Henrion, C.~Prieur, and E.~Tr{\'e}lat, ``Nonlinear optimal control via occupation measures and {LMI}-relaxations,'' \emph{SIAM J. Control Optim.}, vol.~47, no.~4, pp. 1643--1666, 2008.

\bibitem{SCC.Majumdar.Vasudevan.ea2014}
A.~Majumdar, R.~Vasudevan, M.~M. Tobenkin, and R.~Tedrake, ``Convex optimization of nonlinear feedback controllers via occupation measures,'' \emph{Int. J. Robot. Res.}, vol.~33, no.~9, pp. 1209--1230, 2014.

\bibitem{SCC.Korda.Henrion.ea2016}
M.~Korda, D.~Henrion, and C.~N. Jones, ``Controller design and value function approximation for nonlinear dynamical systems,'' \emph{Automatica}, vol.~67, pp. 54 -- 66, 2016.

\bibitem{SCC.Rosenfeld.Kamalapurkar2021}
J.~A. Rosenfeld and R.~Kamalapurkar, ``Dynamic mode decomposition with control {L}iouville operators,'' in \emph{IFAC-PapersOnLine}, vol.~54, no.~9, Jul. 2021, pp. 707--712.

\bibitem{zhdanov2013identity}
F.~Zhdanov and Y.~Kalnishkan, ``An identity for kernel ridge regression,'' \emph{Theoretical Computer Science}, vol. 473, pp. 157--178, 2013.

\bibitem{khalil2002nonlinear}
H.~K. Khalil and J.~W. Grizzle, \emph{Nonlinear systems}.\hskip 1em plus 0.5em minus 0.4em\relax Prentice hall Upper Saddle River, NJ, 2002, vol.~3.

\bibitem{rosenfeld2024dynamic}
J.~A. Rosenfeld and R.~Kamalapurkar, ``Dynamic mode decomposition with control liouville operators,'' \emph{IEEE Transactions on Automatic Control}, 2024.

\bibitem{carmeli2010vector}
C.~Carmeli, E.~De~Vito, A.~Toigo, and V.~Umanit{\'a}, ``Vector valued reproducing kernel hilbert spaces and universality,'' \emph{Analysis and Applications}, vol.~8, no.~01, pp. 19--61, 2010.

\bibitem{SCC.Kimeldorf.Wahba1970}
G.~S. Kimeldorf and G.~Wahba, ``A correspondence between {B}ayesian estimation on stochastic processes and smoothing by splines,'' \emph{Ann. Math. Statist.}, vol.~41, pp. 495--502, 1970.

\bibitem{zhdanov2009competing}
F.~Zhdanov and V.~Vovk, ``Competing with gaussian linear experts,'' ar{X}iv:0910.4683, 2009.

\bibitem{Wahba1990}
G.~Wahba, \emph{Spline Models for Observational Data}.\hskip 1em plus 0.5em minus 0.4em\relax Society for Industrial and Applied Mathematics, 1990.

\bibitem{vovk2013kernel}
V.~Vovk, ``Kernel ridge regression,'' in \emph{Empirical Inference: Festschrift in Honor of Vladimir N. Vapnik}.\hskip 1em plus 0.5em minus 0.4em\relax Springer, 2013, pp. 105--116.

\end{thebibliography}
